\documentclass[12pt]{article}%
\usepackage{amsfonts}
\usepackage{sw20bams}
\usepackage{amsmath}
\usepackage{amssymb}
\usepackage{graphicx}%
\providecommand{\U}[1]{\protect\rule{.1in}{.1in}}
\begin{document}

\title{Uniform Triangles with Equality Constraints}
\author{Steven Finch}
\date{November 26, 2014}
\maketitle

\begin{abstract}
The equality constraint $a+b+c=1$ for random triangle sides corresponds to
breaking a stick in two places. An analog $a^{2}+b^{2}+c^{2}=1$ has a
remarkable feature: the bivariate density for angles coincides with that for
3D\ Gaussian triangles. Interesting complications also arise for $a+b=1$ and
for $a^{2}+b^{2}=1$, with the understanding that the angle $\gamma$ opposite
side $c$ is Uniform$[0,\pi]$. Closed-form expressions for several side moments
remain open.

\end{abstract}

\footnotetext{Copyright \copyright \ 2014 by Steven R. Finch. All rights
reserved.}There is a natural method for generating triangles of unit
perimeter:\ break a stick of length $1$ in two places at random, with the
condition that triangle inequalities are satisfied. \ We denote triangle sides
by $a$, $b$, $c$ and opposite angles by $\alpha$, $\beta$, $\gamma$.
\ Pointers to the literature are found in \cite{Fi1-UnifEqua, IP-UnifEqua};
the bivariate density for two arbitrary sides%
\[
\left\{
\begin{array}
[c]{lll}%
8 &  & \text{if }0<x<1/2\text{, }0<y<1/2\text{ and }x+y>1/2,\\
0 &  & \text{otherwise.}%
\end{array}
\right.
\]
is well-known. \ Proof of the bivariate density for two arbitrary angles
\[
\left\{
\begin{array}
[c]{lll}%
8\,\dfrac{\sin(x)\sin(y)\sin(x+y)}{\left(  \sin(x)+\sin(y)+\sin(x+y)\right)
^{3}} &  & \text{if }0<x<\pi\text{, }0<y<\pi\text{ and }x+y<\pi,\\
0 &  & \text{otherwise}%
\end{array}
\right.
\]
will be given in Section 1. \ The latter is a new result, as far as is known,
although it bears resemblance to formulas in \cite{Fi2-UnifEqua}. \ In
particular, the probability that such a triangle is obtuse is $9-12\ln
(2)\approx0.682$. \ 

Let the lengths of the three pieces (from breaking the stick) instead be
$a^{2}$, $b^{2}$, $c^{2}$. \ Inspiration for this example came from Edelman
\&\ Strang \cite{ES-UnifEqua}. \ We will prove in Section 2 that the bivariate
side density is \ \ \
\[
\left\{
\begin{array}
[c]{lll}%
\dfrac{24\sqrt{3}}{\pi}x\,y &  & \text{if }\left\vert x-y\right\vert
<\sqrt{1-x^{2}-y^{2}}<x+y,\\
0 &  & \text{otherwise}%
\end{array}
\right.
\]
and the bivariate angle density is%
\[
\left\{
\begin{array}
[c]{lll}%
\dfrac{24\sqrt{3}}{\pi}\dfrac{\sin(x)^{2}\sin(y)^{2}\sin(x+y)^{2}}{\left(
\sin(x)^{2}+\sin(y)^{2}+\sin(x+y)^{2}\right)  ^{3}} &  & \text{if }%
0<x<\pi\text{, }0<y<\pi\text{ and }x+y<\pi,\\
0 &  & \text{otherwise.}%
\end{array}
\right.
\]
The probability that such a triangle is obtuse is $1-3\sqrt{3}/(4\pi
)\approx0.586$. \ It is remarkable that angles here are distributed
identically to angles for Gaussian triangles in three-dimensional space
\cite{Fi1-UnifEqua}. \ There is no \textit{a priori} reason to expect such a
coincidence. \ One class of triangles arises synthetically (from breaking a
stick:\ the ambient space doesn't matter)\ while the other class arises
analytically (via a sampling of vertices, \textit{i.e.}, point coordinates in
$\mathbb{R}^{3}$: the ambient space matters).

Let us instead break the stick in just one place at random, giving sides $a$
and $b$. \ Generate independently and uniformly an angle $\gamma$ from the
interval $[0,\pi]$. \ The remaining side $c$ and angles $\alpha$, $\beta$ are
computed via the Law of Cosines. \ We will prove in Section 3 that the
bivariate side density for $a=x$, $c=y$ is \ \ \
\[
\left\{
\begin{array}
[c]{lll}%
\dfrac{2}{\pi}\dfrac{y}{\sqrt{1-y^{2}}\sqrt{4x\left(  1-x\right)  -\left(
1-y^{2}\right)  }} &  & \text{if }\left\vert 2x-1\right\vert <y<1,\\
0 &  & \text{otherwise}%
\end{array}
\right.
\]
and the bivariate angle density for angles $\alpha$, $\beta$ is%
\[
\left\{
\begin{array}
[c]{lll}%
\dfrac{1}{\pi}\dfrac{\sin(x+y)}{\left(  \sin(x)+\sin(y)\right)  ^{2}} &  &
\text{if }0<x<\pi\text{, }0<y<\pi\text{ and }x+y<\pi,\\
0 &  & \text{otherwise.}%
\end{array}
\right.
\]
Thus the side density is complicated while the angle density is simple. \ The
probability that such a triangle is obtuse is $3/2-2/\pi\approx0.863$.

Let the lengths of the two pieces (from breaking the stick) instead be $a^{2}
$, $b^{2}$. \ The angle $\gamma$ is exactly as before. \ We will prove in
Section 4 that the bivariate side density for $a=x$, $c=y$ is%
\[
\left\{
\begin{array}
[c]{lll}%
\dfrac{4}{\pi}\dfrac{x\,y}{\sqrt{4x^{2}\left(  1-x^{2}\right)  -\left(
1-y^{2}\right)  ^{2}}} &  & \text{if }\left\vert x-\sqrt{1-x^{2}}\right\vert
<y<x+\sqrt{1-x^{2}},\\
0 &  & \text{otherwise}%
\end{array}
\right.
\]
and the bivariate angle density for angles $\alpha$, $\beta$ is%
\[
\left\{
\begin{array}
[c]{lll}%
\dfrac{2}{\pi}\dfrac{\sin(x)\sin(y)\sin(x+y)}{\left(  \sin(x)^{2}+\sin
(y)^{2}\right)  ^{2}} &  & \text{if }0<x<\pi\text{, }0<y<\pi\text{ and
}x+y<\pi,\\
0 &  & \text{otherwise.}%
\end{array}
\right.
\]
The probability that such a triangle is obtuse is $3/2-1/\sqrt{2}\approx
0.793$. Another remarkable coincidence occurs here: angles are distributed
identically to angles for pinned Gaussian triangles in two-dimensional space
\cite{Fi2-UnifEqua}.

For variety's sake, select a point $(a,b)$ uniformly on the positive quarter
circle of unit radius, center at $(0,0)$. \ Although $a^{2}+b^{2}=1$ here as
well, sampling from a circle is different from breaking a stick/extracting
square roots. \ The angle $\gamma$ is exactly as before. \ Inspiration for
this example came from Portnoy \cite{Po-UnifEqua}. \ We will prove in Section
5 that the bivariate side density for $a=x$, $c=y$ is \ \ \
\[
\left\{
\begin{array}
[c]{lll}%
\dfrac{4}{\pi^{2}}\dfrac{y}{\sqrt{1-x^{2}}\sqrt{4x^{2}\left(  1-x^{2}\right)
-\left(  1-y^{2}\right)  ^{2}}} &  & \text{if }\left\vert x-\sqrt{1-x^{2}%
}\right\vert <y<x+\sqrt{1-x^{2}},\\
0 &  & \text{otherwise}%
\end{array}
\right.
\]
and the bivariate angle density for angles $\alpha$, $\beta$ is%
\[
\left\{
\begin{array}
[c]{lll}%
\dfrac{2}{\pi^{2}}\dfrac{\sin(x+y)}{\sin(x)^{2}+\sin(y)^{2}} &  & \text{if
}0<x<\pi\text{, }0<y<\pi\text{ and }x+y<\pi,\\
0 &  & \text{otherwise.}%
\end{array}
\right.
\]
Again the side density is complicated while the angle density is simple. \ The
probability that such a triangle is obtuse is $1-\left(  2/\pi^{2}\right)
\ln\left(  1+\sqrt{2}\right)  ^{2}\approx0.842$.

Finally, select a point $(a,b,c)$ uniformly on the positive one-eighth sphere
of unit radius, center at $(0,0,0)$. \ Although $a^{2}+b^{2}+c^{2}=1$ here
like Section 2, sampling from a sphere is different. \ We will prove in
Section 6 that the bivariate side density is
\[
\left\{
\begin{array}
[c]{lll}%
\dfrac{C}{\sqrt{x^{2}+y^{2}}\sqrt{1-x^{2}-y^{2}}} &  & \text{if }\left\vert
x-y\right\vert <\sqrt{1-x^{2}-y^{2}}<x+y,\\
0 &  & \text{otherwise}%
\end{array}
\right.
\]
and the bivariate angle density is%
\[
\left\{
\begin{array}
[c]{lll}%
\dfrac{C\sin(x)\sin(y)\sin(x+y)}{\left(  \sin(x)^{2}+\sin(y)^{2}+\sin
(x+y)^{2}\right)  \sqrt{\sin(x)^{2}+\sin(y)^{2}}} &  & \text{if }%
0<x<\pi\text{, }0<y<\pi\text{ and }x+y<\pi,\\
0 &  & \text{otherwise}%
\end{array}
\right.
\]
where $C$ is a known constant (written in terms of dilogarithm function
values). The probability that such a triangle is obtuse is approximately
$0.659$. \ This is the first of several quantities appearing here for which
explicit formulation is unavailable. One quantity, given as approximately
$0.958$ in Section 5, is especially important to understand more fully.
Insight and help toward unraveling these constants would be appreciated.

\section{Constraint $a+b+c=1$}

If we break a line segment in two places at random, the three pieces can be
configured as a triangle with probability $1/4$. \ Reason:\ the subdomain
$\{(x,y):0<x<1/2$, $0<y<1/2$, $x+y>1/2\}$ occupies one-fourth the area of
domain $\{(x,y):0<x<1$, $0<y<1$, $x+y<1\}$, and the triangle inequalities%
\[%
\begin{array}
[c]{ccccc}%
1-a-b=c<a+b, &  & a<b+c=1-a, &  & b<a+c=1-b
\end{array}
\]
become $a+b>1/2$, $a<1/2$ and $b<1/2$. \ The Law of Sines gives%
\[
b\sin(\alpha)-a\sin(\beta)=0
\]
and the Law of Cosines gives%
\begin{align*}
b\cos(\alpha)+a\cos(\beta)  &  =\frac{-a^{2}+b^{2}+(1-a-b)^{2}}{2(1-a-b)}%
+\frac{a^{2}-b^{2}+(1-a-b)^{2}}{2(1-a-b)}\\
&  =1-a-b.
\end{align*}
Solving for $a$, $b$ yields%
\[%
\begin{array}
[c]{ccc}%
a=\dfrac{\sin(\alpha)}{\sin(\alpha)+\sin(\beta)+\sin(\alpha+\beta)}, &  &
b=\dfrac{\sin(\beta)}{\sin(\alpha)+\sin(\beta)+\sin(\alpha+\beta)}%
\end{array}
\]
and thus%
\[
1-a-b=\dfrac{\sin(\alpha+\beta)}{\sin(\alpha)+\sin(\beta)+\sin(\alpha+\beta
)}.
\]
Now, the map $(a,b)\mapsto(\alpha,\beta)$ defined via the Law of Cosines has
Jacobian determinant%
\[
\left\vert J\right\vert =\frac{1}{a\,b(1-a-b)}
\]
hence the desired bivariate density for angles is%
\[
\frac{8}{\left\vert J\right\vert }=8\,\dfrac{\sin(\alpha)\sin(\beta
)\sin(\alpha+\beta)}{\left(  \sin(\alpha)+\sin(\beta)+\sin(\alpha
+\beta)\right)  ^{3}}.
\]
We have univariate densities%
\[
\left\{
\begin{array}
[c]{lll}%
8\,a &  & \text{if }0<a<1/2\text{,}\\
0 &  & \text{otherwise;}%
\end{array}
\right.
\]%
\[
\left\{
\begin{array}
[c]{lll}%
-8\,\dfrac{(3-\cos(\alpha))\sin(\alpha)}{(1+\cos(\alpha))^{3}}\ln\left(
\sin\left(  \dfrac{\alpha}{2}\right)  \right)  -8\,\dfrac{\sin(\alpha
)}{(1+\cos(\alpha))^{2}} &  & \text{if }0<\alpha<\pi\text{,}\\
0 &  & \text{otherwise}%
\end{array}
\right.
\]
and moments%
\[%
\begin{array}
[c]{ccccc}%
\operatorname*{E}(a)=1/3, &  & \operatorname*{E}(a^{2})=1/8, &  &
\operatorname*{E}(a\,b)=5/48;
\end{array}
\]%
\[%
\begin{array}
[c]{ccccc}%
\operatorname*{E}(\alpha)=\pi/3, &  & \operatorname*{E}(\alpha^{2}%
)=8/3-\pi^{2}/9, &  & \operatorname*{E}(\alpha\,\beta)=-4/3+2\pi^{2}/9.
\end{array}
\]%
\begin{figure}[ptb]%
\centering
\includegraphics[
height=2.9637in,
width=6.2872in
]%
{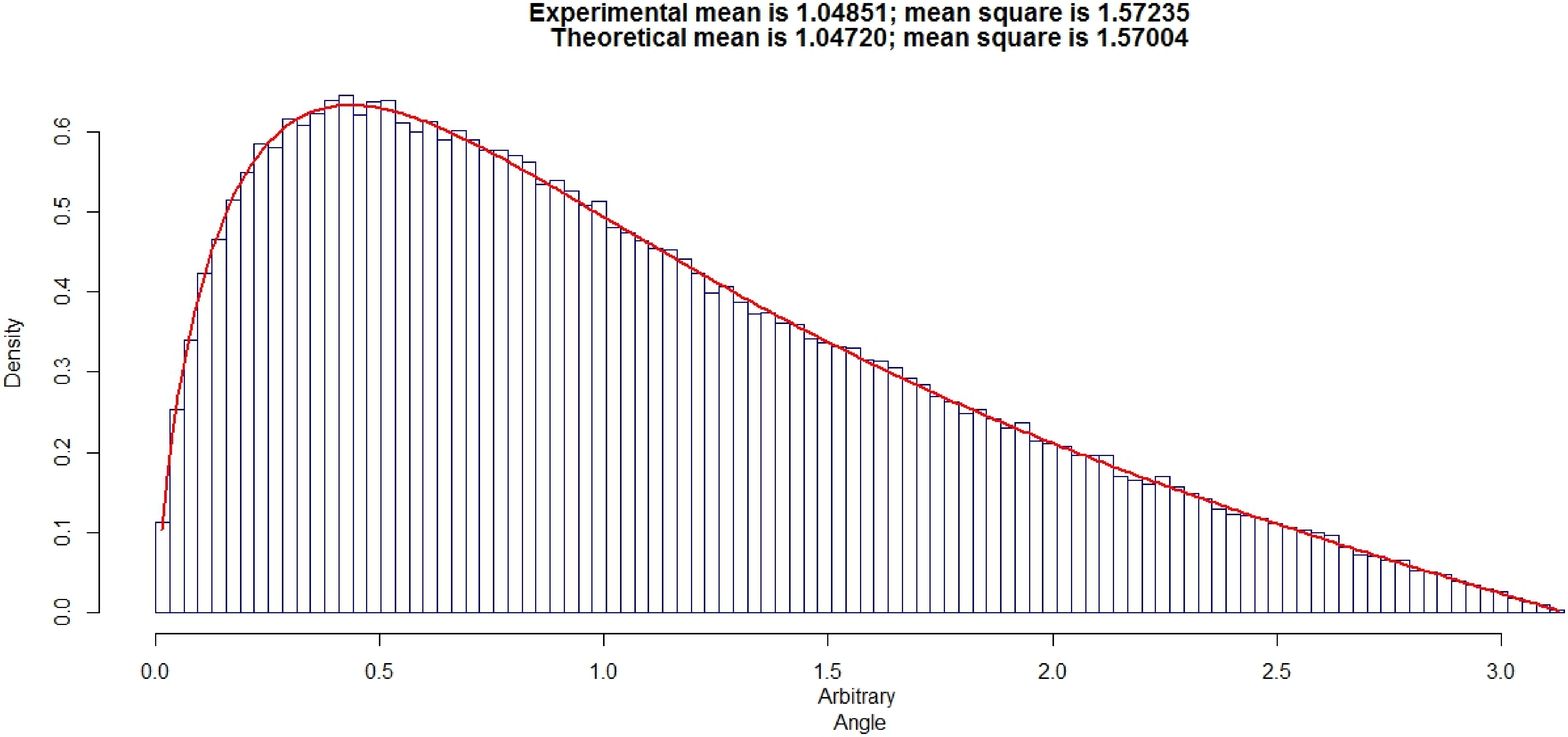}%
\caption{Density function for arbitrary angle in Section 1.}%
\end{figure}

\section{Constraint $a^{2}+b^{2}+c^{2}=1$}

Here $a^{2}$, $b^{2}$ are uniform, not $a$, $b$. \ The triangle inequalities%
\[
\sqrt{1-a^{2}-b^{2}}=c<a+b,
\]%
\[
a<b+c=b+\sqrt{1-a^{2}-b^{2}},
\]%
\[
b<a+c=a+\sqrt{1-a^{2}-b^{2}}
\]
are simultaneously met with probability $\sqrt{3}\pi/9\approx0.604$
\cite{ES-UnifEqua} because the subdomain%
\[
\left\{  (x,y):\left\vert x-y\right\vert <\sqrt{1-x^{2}-y^{2}}<x+y\right\}
\]
has area $\sqrt{3}\pi/18$ whereas the domain $\{(x,y):0<x<1$, $0<y<1$,
$x+y<1\}$ has area $1/2$. The map $(u,v)\mapsto(\sqrt{u},\sqrt{v})$ has
Jacobian determinant $1/(4\sqrt{u\,v})$, providing the form $a\,b$ of the
bivariate side density. As before, the Law of Sines gives%
\[
b\sin(\alpha)-a\sin(\beta)=0
\]
and the Law of Cosines gives%
\begin{align*}
b\cos(\alpha)+a\cos(\beta)  &  =\frac{-a^{2}+b^{2}+(1-a^{2}-b^{2})}%
{2\sqrt{1-a^{2}-b^{2}}}+\frac{a^{2}-b^{2}+(1-a^{2}-b^{2})}{2\sqrt
{1-a^{2}-b^{2}}}\\
&  =\sqrt{1-a^{2}-b^{2}}.
\end{align*}
Solving for $a$, $b$ yields%
\[%
\begin{array}
[c]{ccc}%
a=\dfrac{\sin(\alpha)}{\sqrt{\sin(\alpha)^{2}+\sin(\beta)^{2}+\sin
(\alpha+\beta)^{2}}}, &  & b=\dfrac{\sin(\beta)}{\sqrt{\sin(\alpha)^{2}%
+\sin(\beta)^{2}+\sin(\alpha+\beta)^{2}}}%
\end{array}
\]
and thus%
\[
\sqrt{1-a^{2}-b^{2}}=\dfrac{\sin(\alpha+\beta)}{\sqrt{\sin(\alpha)^{2}%
+\sin(\beta)^{2}+\sin(\alpha+\beta)^{2}}}.
\]
Now, the map $(a,b)\mapsto(\alpha,\beta)$ defined via the Law of Cosines has
Jacobian determinant%
\[
\left\vert J\right\vert =\frac{1}{a\,b(1-a^{2}-b^{2})}
\]
hence the desired bivariate density for angles is%
\[
\dfrac{24\sqrt{3}}{\pi}\frac{a\,b}{\left\vert J\right\vert }=\dfrac{24\sqrt
{3}}{\pi}\dfrac{\sin(\alpha)^{2}\sin(\beta)^{2}\sin(\alpha+\beta)^{2}}{\left(
\sin(\alpha)^{2}+\sin(\beta)^{2}+\sin(\alpha+\beta)^{2}\right)  ^{3}}.
\]
We have univariate densities%
\[
\left\{
\begin{array}
[c]{lll}%
\dfrac{12\sqrt{3}}{\pi}a^{2}\sqrt{2-3a^{2}} &  & \text{if }0<a<\dfrac{\sqrt
{6}}{3}\text{,}\\
0 &  & \text{otherwise;}%
\end{array}
\right.
\]%
\[
\left\{
\begin{array}
[c]{lll}%
\dfrac{6\sqrt{3}}{\pi}\dfrac{\left(  2+\cos(\alpha)^{2}\right)  \sin(\alpha
)}{\left(  4-\cos(\alpha)^{2}\right)  ^{5/2}}\left(  \dfrac{\pi}{2}%
+\arcsin\left(  \dfrac{\cos(\alpha)}{2}\right)  \right)  +\dfrac{9\sqrt{3}%
}{\pi}\dfrac{\cos(\alpha)\sin(\alpha)}{\left(  4-\cos(\alpha)^{2}\right)
^{2}} &  & \text{if }0<\alpha<\pi\text{,}\\
0 &  & \text{otherwise}%
\end{array}
\right.
\]
and moments%
\[%
\begin{array}
[c]{ccccc}%
\operatorname*{E}(a)=\dfrac{32\sqrt{6}}{45\pi}, &  & \operatorname*{E}%
(a^{2})=\dfrac{1}{3}, &  & \operatorname*{E}(a\,b)=\dfrac{9+\sqrt{3}\pi
}{9\sqrt{3}\pi};
\end{array}
\]%
\[%
\begin{array}
[c]{ccccc}%
\operatorname*{E}(\alpha)=\dfrac{\pi}{3}, &  & \operatorname*{E}(\alpha
^{2})=\dfrac{1}{3}\left(  \pi-\sqrt{3}\right)  \pi, &  & \operatorname*{E}%
(\alpha\,\beta)=\dfrac{\sqrt{3}\pi}{6}.
\end{array}
\]
We are familiar with the angle density \cite{Fi2-UnifEqua}; the side density,
however, is new.

\section{Constraint $a+b=1$ and Uniform $\gamma$}

The map $(a,\gamma)\mapsto(a,c)$ defined by%
\[
c^{2}=a^{2}+(1-a)^{2}-2a(1-a)\cos(\gamma)
\]
has Jacobian determinant%
\begin{align*}
\left\vert J\right\vert  &  =\frac{\partial c}{\partial\gamma}=\frac{1}%
{2c}\cdot2a(1-a)\sin(\gamma)\\
&  =\frac{a(1-a)}{c}\sqrt{1-\cos(\gamma)^{2}}\\
&  =\frac{a(1-a)}{c}\sqrt{1-\left[  \frac{a^{2}+(1-a)^{2}-c^{2}}%
{2a(1-a)}\right]  ^{2}}\\
&  =\frac{a(1-a)}{c}\sqrt{\frac{4a^{2}(1-a)^{2}-\left[  a^{2}+(1-a)^{2}%
-c^{2}\right]  ^{2}}{4a^{2}(1-a)^{2}}}\\
&  =\frac{1}{2c}\sqrt{4a^{2}(1-a)^{2}-\left[  2a^{2}-2a+1-c^{2}\right]  ^{2}%
}\\
&  =\frac{1}{2c}\sqrt{4a^{2}(1-a)^{2}-\left[  -2a(1-a)+\left(  1-c^{2}\right)
\right]  ^{2}}\\
&  =\frac{1}{2c}\sqrt{4a(1-a)\left(  1-c^{2}\right)  -\left(  1-c^{2}\right)
^{2}}\\
&  =\frac{1}{2c}\sqrt{1-c^{2}}\sqrt{4a(1-a)-\left(  1-c^{2}\right)  };
\end{align*}
therefore the bivariate side density is%
\[
\frac{1}{\pi}\frac{1}{\left\vert J\right\vert }=\frac{2}{\pi}\frac{c}%
{\sqrt{1-c^{2}}\sqrt{4a(1-a)-\left(  1-c^{2}\right)  }}.
\]
The Law of Sines gives%
\[
c\sin(\alpha)-a\sin(\alpha+\beta)=0
\]
and the Law of Cosines gives%
\begin{align*}
\left(  1-a\right)  \cos(\alpha)+a\cos(\beta)  &  =b\cos(\alpha)+a\cos
(\beta)\\
&  =\frac{-a^{2}+b^{2}+c^{2}}{2c}+\frac{a^{2}-b^{2}+c^{2}}{2c}=c.
\end{align*}
Solving for $a$, $c$ yields%
\[%
\begin{array}
[c]{ccc}%
a=\dfrac{\sin(\alpha)}{\sin(\alpha)+\sin(\beta)}, &  & c=\dfrac{\sin
(\alpha+\beta)}{\sin(\alpha)+\sin(\beta)}%
\end{array}
\]
and thus%
\[
b=1-a=\dfrac{\sin(\beta)}{\sin(\alpha)+\sin(\beta)}.
\]
Now, the map $(a,c)\mapsto(\alpha,\beta)$ defined via the Law of Cosines has
Jacobian determinant%
\[
\left\vert I\right\vert =\frac{1}{a\,b\,c}
\]
hence the desired bivariate density for angles is%
\begin{align*}
\dfrac{2}{\pi}\dfrac{c}{\sqrt{1-c^{2}}\sqrt{4a\left(  1-a\right)  -\left(
1-c^{2}\right)  }}\frac{1}{\left\vert I\right\vert }  &  =\dfrac{1}{\pi}%
\dfrac{\sin(\alpha)+\sin(\beta)}{\sin(\alpha)\sin(\beta)}\dfrac{\sin
(\alpha)\sin(\beta)\sin(\alpha+\beta)}{\left(  \sin(\alpha)+\sin
(\beta)\right)  ^{3}}\\
&  =\dfrac{1}{\pi}\dfrac{\sin(\alpha+\beta)}{\left(  \sin(\alpha)+\sin
(\beta)\right)  ^{2}}.
\end{align*}
We have univariate densities%
\[%
\begin{array}
[c]{ccc}%
\left\{
\begin{array}
[c]{lll}%
1 &  & \text{if }0<a<1\text{,}\\
0 &  & \text{otherwise;}%
\end{array}
\right.  &  & \left\{
\begin{array}
[c]{lll}%
\dfrac{c}{\sqrt{1-c^{2}}} &  & \text{if }0<c<1\text{,}\\
0 &  & \text{otherwise;}%
\end{array}
\right.
\end{array}
\]%
\[
\left\{
\begin{array}
[c]{lll}%
-\dfrac{1}{\pi}\dfrac{1}{\cos(\alpha)^{2}}\left[  2\ln\left(  \sin\left(
\dfrac{\alpha}{2}\right)  \right)  +\ln(2)\right]  -\dfrac{1}{\pi}\dfrac
{1}{\cos(\alpha)} &  & \text{if }0<\alpha<\pi\text{,}\\
0 &  & \text{otherwise}%
\end{array}
\right.
\]
and moments%
\[%
\begin{array}
[c]{ccccccccc}%
\operatorname*{E}(a)=\dfrac{1}{2}, &  & \operatorname*{E}(a^{2})=\dfrac{1}%
{3}, &  & \operatorname*{E}(c)=\dfrac{\pi}{4}, &  & \operatorname*{E}%
(c^{2})=\dfrac{2}{3}, &  & \operatorname*{E}(a\,c)=\dfrac{\pi}{8};
\end{array}
\]%
\[%
\begin{array}
[c]{ccccc}%
\operatorname*{E}(\alpha)=\dfrac{\pi}{4}, &  & \operatorname*{E}(\alpha
^{2})=1.3029200473..., &  & \operatorname*{E}(\alpha\,\beta)=0.3420140195....
\end{array}
\]

\section{Constraint $a^{2}+b^{2}=1$ and Uniform $\gamma$}

The map $(a,\gamma)\mapsto(a,c)$ defined by%
\[
c^{2}=a^{2}+(1-a^{2})-2a\sqrt{1-a^{2}}\cos(\gamma)
\]
has Jacobian determinant%
\begin{align*}
\left\vert J\right\vert  &  =\frac{\partial c}{\partial\gamma}=\frac{1}%
{2c}\cdot2a\sqrt{1-a^{2}}\sin(\gamma)\\
&  =\frac{a\sqrt{1-a^{2}}}{c}\sqrt{1-\cos(\gamma)^{2}}\\
&  =\frac{a\sqrt{1-a^{2}}}{c}\sqrt{1-\left[  \frac{a^{2}+(1-a^{2})-c^{2}%
}{2a\sqrt{1-a^{2}}}\right]  ^{2}}\\
&  =\frac{a\sqrt{1-a^{2}}}{c}\sqrt{\frac{4a^{2}(1-a^{2})-\left(
1-c^{2}\right)  ^{2}}{4a^{2}(1-a^{2})}}\\
&  =\frac{1}{2c}\sqrt{4a^{2}(1-a^{2})-\left(  1-c^{2}\right)  ^{2}};
\end{align*}
therefore the bivariate side density is%
\[
\frac{1}{\pi}\frac{4a}{\left\vert J\right\vert }=\dfrac{4}{\pi}\dfrac
{a\,c}{\sqrt{4a^{2}\left(  1-a^{2}\right)  -\left(  1-c^{2}\right)  ^{2}}}.
\]
The additional factor $a$ comes because $u\mapsto\sqrt{u}$ has derivative
$1/(2\sqrt{u})$. As before, the Law of Sines gives%
\[
c\sin(\alpha)-a\sin(\alpha+\beta)=0
\]
and the Law of Cosines gives%
\begin{align*}
\sqrt{1-a^{2}}\cos(\alpha)+a\cos(\beta)  &  =b\cos(\alpha)+a\cos(\beta)\\
&  =\frac{-a^{2}+b^{2}+c^{2}}{2c}+\frac{a^{2}-b^{2}+c^{2}}{2c}=c.
\end{align*}
Solving for $a$, $c$ yields%
\[%
\begin{array}
[c]{ccc}%
a=\dfrac{\sin(\alpha)}{\sqrt{\sin(\alpha)^{2}+\sin(\beta)^{2}}}, &  &
c=\dfrac{\sin(\alpha+\beta)}{\sqrt{\sin(\alpha)^{2}+\sin(\beta)^{2}}}%
\end{array}
\]
and thus%
\[
b^{2}=1-a^{2}=\dfrac{\sin(\beta)^{2}}{\sin(\alpha)^{2}+\sin(\beta)^{2}}.
\]
Now, the map $(a,c)\mapsto(\alpha,\beta)$ defined via the Law of Cosines has
Jacobian determinant%
\[
\left\vert I\right\vert =\frac{1}{a\,b^{2}c}
\]
hence the desired bivariate density for angles is%
\begin{align*}
\dfrac{4}{\pi}\dfrac{a\,c}{\sqrt{4a^{2}\left(  1-a^{2}\right)  -\left(
1-c^{2}\right)  ^{2}}}\frac{1}{\left\vert I\right\vert }  &  =\dfrac{2}{\pi
}\frac{1}{\sin(\beta)}\dfrac{\sin(\alpha)\sin(\beta)^{2}\sin(\alpha+\beta
)}{\left(  \sin(\alpha)^{2}+\sin(\beta)^{2}\right)  ^{2}}\\
&  =\dfrac{2}{\pi}\dfrac{\sin(\alpha)\sin(\beta)\sin(\alpha+\beta)}{\left(
\sin(\alpha)^{2}+\sin(\beta)^{2}\right)  ^{2}}.
\end{align*}
We have univariate densities%
\[%
\begin{array}
[c]{ccc}%
\left\{
\begin{array}
[c]{lll}%
2a &  & \text{if }0<a<1\text{,}\\
0 &  & \text{otherwise;}%
\end{array}
\right.  &  & \left\{
\begin{array}
[c]{lll}%
c &  & \text{if }0<c<\sqrt{2}\text{,}\\
0 &  & \text{otherwise;}%
\end{array}
\right.
\end{array}
\]%
\[
\left\{
\begin{array}
[c]{lll}%
\dfrac{1}{\pi}\dfrac{\cos(\alpha)}{\left(  2-\cos(\alpha)^{2}\right)  ^{3/2}%
}\left(  \dfrac{\pi}{2}+\arcsin\left(  \dfrac{\cos(\alpha)}{\sqrt{2}}\right)
\right)  +\dfrac{1}{\pi}\dfrac{1}{2-\cos(\alpha)^{2}} &  & \text{if }%
0<\alpha<\pi\text{,}\\
0 &  & \text{otherwise}%
\end{array}
\right.
\]
and moments%
\[%
\begin{array}
[c]{ccccccc}%
\operatorname*{E}(a)=\dfrac{2}{3}, &  & \operatorname*{E}(a^{2})=\dfrac{1}%
{2}, &  & \operatorname*{E}(c)=\dfrac{2\sqrt{2}}{3}, &  & \operatorname*{E}%
(c^{2})=1,
\end{array}
\]%
\[
\operatorname*{E}(a\,c)=\frac{\sqrt{2}}{\pi}\,\,%
{\displaystyle\int\limits_{0}^{\sqrt{2}}}
t^{2}\sqrt{1+t\sqrt{2-t^{2}}}\,E\left(  \sqrt{\frac{2t\sqrt{2-t^{2}}}%
{1+t\sqrt{2-t^{2}}}}\right)  dt=0.6272922529...
\]
where
\[
E(\xi)=%
{\displaystyle\int\limits_{0}^{\pi/2}}
\sqrt{1-\xi^{2}\sin(\theta)^{2}}\,d\theta=%
{\displaystyle\int\limits_{0}^{1}}
\sqrt{\dfrac{1-\xi^{2}t^{2}}{1-t^{2}}}\,dt
\]
is the complete elliptic integral of the second kind;
\[%
\begin{array}
[c]{ccccc}%
\operatorname*{E}(\alpha)=\dfrac{\pi}{4}, &  & \operatorname*{E}(\alpha
^{2})=\dfrac{5}{48}\pi^{2}+\dfrac{1}{4}\ln\left(  2\right)  ^{2}, &  &
\operatorname*{E}(\alpha\,\beta)=\dfrac{1}{16}\pi^{2}-\dfrac{1}{4}\ln\left(
2\right)  ^{2}.
\end{array}
\]
We are familiar with the angle density \cite{Fi2-UnifEqua}; the side density,
however, is new.

\section{Constraint $a^{2}+b^{2}=1$ Revisited}

Let $\theta$ be uniformly distributed on the interval $[0,\pi/2]$. The map
$(\theta,\gamma)\mapsto(a,c)$ defined by%
\[
\left\{
\begin{array}
[c]{l}%
a=\cos(\theta),\\
c^{2}=a^{2}+(1-a^{2})-2a\sqrt{1-a^{2}}\cos(\gamma)
\end{array}
\right.
\]
has Jacobian determinant%
\begin{align*}
\left\vert J\right\vert  &  =\frac{\partial a}{\partial\theta}\frac{\partial
c}{\partial\gamma}=\sin(\theta)\cdot\frac{1}{2c}\cdot2a\sqrt{1-a^{2}}%
\sin(\gamma)\\
&  =\frac{a\left(  1-a^{2}\right)  }{c}\sqrt{1-\cos(\gamma)^{2}}\\
&  =\frac{a\left(  1-a^{2}\right)  }{c}\sqrt{1-\left[  \frac{a^{2}%
+(1-a^{2})-c^{2}}{2a\sqrt{1-a^{2}}}\right]  ^{2}}\\
&  =\frac{a\left(  1-a^{2}\right)  }{c}\sqrt{\frac{4a^{2}(1-a^{2})-\left(
1-c^{2}\right)  ^{2}}{4a^{2}(1-a^{2})}}\\
&  =\frac{1}{2c}\sqrt{1-a^{2}}\sqrt{4a^{2}(1-a^{2})-\left(  1-c^{2}\right)
^{2}};
\end{align*}
therefore the bivariate side density is%
\[
\frac{2}{\pi^{2}}\frac{1}{\left\vert J\right\vert }=\dfrac{4}{\pi^{2}}%
\dfrac{c}{\sqrt{1-a^{2}}\sqrt{4a^{2}\left(  1-a^{2}\right)  -\left(
1-c^{2}\right)  ^{2}}}.
\]
\ The same expressions for $a$, $b$, $c$, $\left\vert I\right\vert $ apply as
in Section 4, hence the desired bivariate density for angles is%
\begin{align*}
\dfrac{4}{\pi^{2}}\dfrac{c}{\sqrt{1-a^{2}}\sqrt{4a^{2}\left(  1-a^{2}\right)
-\left(  1-c^{2}\right)  ^{2}}}\frac{1}{\left\vert I\right\vert }  &
=\dfrac{2}{\pi^{2}}\frac{\sin(\alpha)^{2}+\sin(\beta)^{2}}{\sin(\alpha
)\sin(\beta)^{2}}\dfrac{\sin(\alpha)\sin(\beta)^{2}\sin(\alpha+\beta)}{\left(
\sin(\alpha)^{2}+\sin(\beta)^{2}\right)  ^{2}}\\
&  =\dfrac{2}{\pi^{2}}\dfrac{\sin(\alpha+\beta)}{\sin(\alpha)^{2}+\sin
(\beta)^{2}}.
\end{align*}
We have univariate densities%
\[%
\begin{array}
[c]{ccc}%
\left\{
\begin{array}
[c]{lll}%
\dfrac{2}{\pi}\dfrac{1}{\sqrt{1-a^{2}}} &  & \text{if }0<a<1\text{,}\\
0 &  & \text{otherwise;}%
\end{array}
\right.  &  & \left\{
\begin{array}
[c]{lll}%
\dfrac{4c}{\pi^{2}}K\left(  c\sqrt{2-c^{2}}\right)  &  & \text{if }%
0<c<\sqrt{2}\text{,}\\
0 &  & \text{otherwise}%
\end{array}
\right.
\end{array}
\]
where
\[
K(\xi)=%
{\displaystyle\int\limits_{0}^{\pi/2}}
\dfrac{1}{\sqrt{1-\xi^{2}\sin(\theta)^{2}}}\,d\theta=%
{\displaystyle\int\limits_{0}^{1}}
\dfrac{1}{\sqrt{(1-t^{2})(1-\xi^{2}t^{2})}}\,dt
\]
is the complete elliptic integral of the first kind;%
\[
\left\{
\begin{array}
[c]{lll}%
\dfrac{1}{2\pi}+\dfrac{1}{\pi^{2}}\dfrac{\cos(\alpha)}{\sqrt{2-\cos
(\alpha)^{2}}}\ln\left(  \dfrac{2-\cos(\alpha)+\sqrt{2-\cos(\alpha)^{2}}%
}{2-\cos(\alpha)-\sqrt{2-\cos(\alpha)^{2}}}\right)  &  & \text{if }%
0<\alpha<\pi\text{,}\\
0 &  & \text{otherwise}%
\end{array}
\right.
\]
and moments%
\[%
\begin{array}
[c]{ccccccc}%
\operatorname*{E}(a)=\dfrac{2}{\pi}, &  & \operatorname*{E}(a^{2})=\dfrac
{1}{2}, &  & \operatorname*{E}(c)=0.9580913986..., &  & \operatorname*{E}%
(c^{2})=1,
\end{array}
\]%
\[
\operatorname*{E}(a\,c)=\frac{2\sqrt{2}}{\pi^{2}}\,\,%
{\displaystyle\int\limits_{0}^{\sqrt{2}}}
\frac{t^{2}}{\sqrt{1+t\sqrt{2-t^{2}}}}\,K\left(  \sqrt{\frac{2t\sqrt{2-t^{2}}%
}{1+t\sqrt{2-t^{2}}}}\right)  dt=0.6080033617...,
\]%
\[%
\begin{array}
[c]{ccccc}%
\operatorname*{E}(\alpha)=\dfrac{\pi}{4}, &  & \operatorname*{E}(\alpha
^{2})=1.2565739217..., &  & \operatorname*{E}(\alpha\,\beta)=0.3883601451....
\end{array}
\]
Evaluating the mean of side $c$ in closed-form remains tantalizingly open.
\ In addition to its definition:%
\[
\operatorname*{E}(c)=\dfrac{4}{\pi^{2}}%
{\displaystyle\int\limits_{0}^{\sqrt{2}}}
\,t^{2}K\left(  t\sqrt{2-t^{2}}\right)  dt
\]
we have the following representation:%
\[
\operatorname*{E}(c)=\dfrac{2}{\pi^{2}}%
{\displaystyle\int\limits_{0}^{1}}
s\frac{\sqrt{1-\sqrt{1-s^{2}}}+\sqrt{1+\sqrt{1-s^{2}}}}{\sqrt{1-s^{2}}%
}K(s)\,ds
\]
which unfortunately does not appear in \cite{Gs-UnifEqua}.%
\begin{figure}[ptb]%
\centering
\includegraphics[
trim=0.000000in -0.077642in 0.000000in 0.077642in,
height=2.821in,
width=6.4022in
]%
{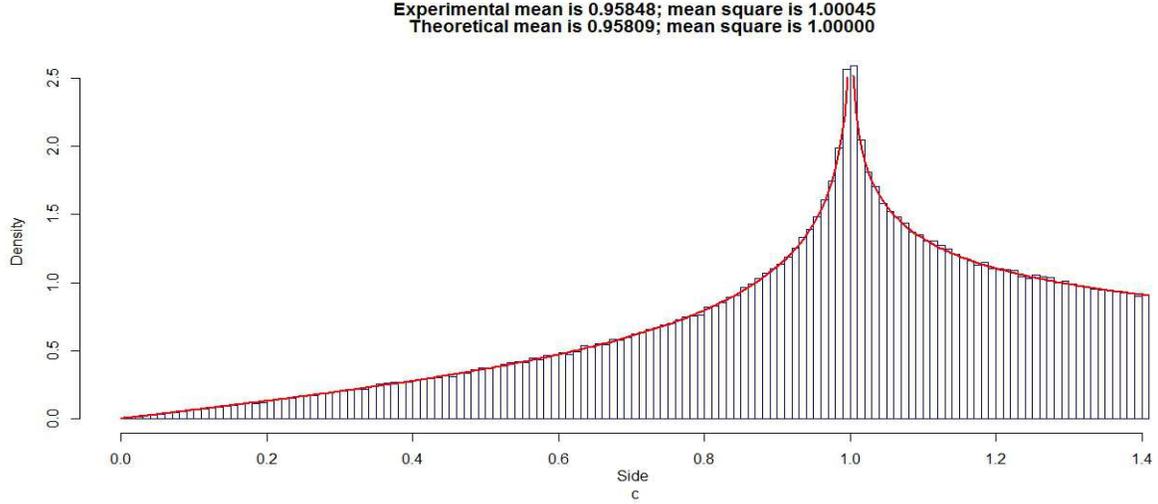}%
\caption{Density function for side $c$ in Section 5.}%
\end{figure}

\section{Constraint $a^{2}+b^{2}+c^{2}=1$ Revisited}

Let $\varphi$, $\psi$ be independently and uniformly distributed on the
interval $[0,\pi/2]$. The map $(\varphi,\psi)\mapsto(a,b)$ defined by%
\[
\left\{
\begin{array}
[c]{l}%
a=\sin(\varphi)\cos(\psi),\\
b=\sin(\varphi)\sin(\psi)
\end{array}
\right.
\]
has Jacobian determinant%
\[
\left\vert J\right\vert =\cos(\varphi)\sin(\varphi)=\sqrt{a^{2}+b^{2}}%
\sqrt{1-a^{2}-b^{2}}
\]
therefore the bivariate side density is%
\[
\frac{C}{\left\vert J\right\vert }=\dfrac{C}{\sqrt{a^{2}+b^{2}}\sqrt
{1-a^{2}-b^{2}}}
\]
and the normalizing constant satisfies%
\begin{align*}
\frac{1}{C}  &  =2%
{\displaystyle\int\limits_{0}^{1}}
\frac{\arctan\left(  x\sqrt{1+x^{2}}\right)  }{1+x^{2}}dx\\
&  =\frac{\pi}{2}\arctan\left(  \sqrt{2}\right)  -2%
{\displaystyle\int\limits_{0}^{1}}
\frac{1+2x^{2}}{\sqrt{1+x^{2}}\left(  1+x^{2}+x^{4}\right)  }dx\\
&  =-\frac{\pi^{2}}{24}-\frac{1}{2}i\,\pi\ln\left(  2-\sqrt{3}\right)
-\operatorname*{Li}\nolimits_{2}\left(  \sqrt{2-\sqrt{3}}\right)  +\\
&  \operatorname*{Li}\nolimits_{2}\left(  -\sqrt{2-\sqrt{3}}\right)
-\operatorname*{Li}\nolimits_{2}\left(  -\sqrt{2+\sqrt{3}}\right)
+\operatorname*{Li}\nolimits_{2}\left(  \sqrt{2+\sqrt{3}}\right)  +\\
&  \operatorname*{Li}\nolimits_{2}\left(  \frac{1-i}{1+\sqrt{3}}\right)
+\operatorname*{Li}\nolimits_{2}\left(  \frac{1+i}{1+\sqrt{3}}\right)
+\operatorname*{Li}\nolimits_{2}\left(  \frac{1+\sqrt{3}}{-1+i}\right)
+\operatorname*{Li}\nolimits_{2}\left(  \frac{1+\sqrt{3}}{-1-i}\right) \\
&  =0.6947951075...
\end{align*}
where%
\[
\operatorname*{Li}\nolimits_{2}(x)=%
{\displaystyle\sum\limits_{k=1}^{\infty}}
\dfrac{x^{k}}{k^{2}}=-%
{\displaystyle\int\limits_{0}^{x}}
\dfrac{\ln(1-t)}{t}dt
\]
is the dilogarithm function. \ The same expressions for $a$, $b$, $c$,
$\left\vert I\right\vert $ apply as in Section 2, hence the desired bivariate
density for angles is%
\[
\dfrac{C\,a\,b\sqrt{1-a^{2}-b^{2}}}{\sqrt{a^{2}+b^{2}}}=\dfrac{C\sin
(\alpha)\sin(\beta)\sin(\alpha+\beta)}{\left(  \sin(\alpha)^{2}+\sin
(\beta)^{2}+\sin(\alpha+\beta)^{2}\right)  \sqrt{\sin(\alpha)^{2}+\sin
(\beta)^{2}}}.
\]
Let%
\[
\omega(a)=\arcsin\left(  \frac{a+\sqrt{2-3a^{2}}}{\sqrt{1-a^{2}}\sqrt
{4a^{2}+\left(  a+\sqrt{2-3a^{2}}\right)  ^{2}}}\right)
\]
then the univariate side density is%
\[
\left\{
\begin{array}
[c]{lll}%
-C\,F\left(  \omega(-a),\sqrt{1-a^{2}}\right)  +C\,F\left(  \omega
(a),\sqrt{1-a^{2}}\right)  &  & \text{if }0<a<\sqrt{1/2}\text{,}\\
-C\,F\left(  -\omega(-a),\sqrt{1-a^{2}}\right)  +C\,F\left(  \omega
(a),\sqrt{1-a^{2}}\right)  &  & \text{if }\sqrt{1/2}<a<\sqrt{2/3}\text{,}\\
0 &  & \text{otherwise}%
\end{array}
\right.
\]
where
\[
F(\omega,\xi)=%
{\displaystyle\int\limits_{0}^{\omega}}
\dfrac{1}{\sqrt{1-\xi^{2}\sin(\theta)^{2}}}\,d\theta=%
{\displaystyle\int\limits_{0}^{\sin(\omega)}}
\dfrac{1}{\sqrt{(1-t^{2})(1-\xi^{2}t^{2})}}\,dt
\]
is the incomplete elliptic integral of the first kind. \ We have not attempted
to evaluate the univariate angle density; numerical calculations lead to
moments:%
\[%
\begin{array}
[c]{ccccc}%
\operatorname*{E}(a)=0.5361308550..., &  & \operatorname*{E}(a^{2}%
)=0.3209403207..., &  & \operatorname*{E}(a\,b)=0.2707436816...,
\end{array}
\]%
\[%
\begin{array}
[c]{ccccc}%
\operatorname*{E}(\alpha)=1.0018939715..., &  & \operatorname*{E}(\alpha
^{2})=1.4360872743..., &  & \operatorname*{E}(\alpha\,\beta)=0.8093206054...
\end{array}
\]
none of which are immediately recognizable.

For the sake of thoroughness (but without proof), the pieces $a$, $b$, $c$ can
be configured as a triangle with probability%
\begin{align*}
\Delta &  =\frac{8}{\pi^{2}}%
{\displaystyle\int\limits_{\lambda}^{\pi/4}}
\left[  \frac{\pi}{4}-2\arctan\left(  \frac{\sin(\varphi)-\sqrt{-\cos
(\varphi)^{2}+2\sin(\varphi)^{2}}}{\cos(\varphi)+\sin(\varphi)}\right)
\right]  d\varphi+\\
&  \frac{8}{\pi^{2}}%
{\displaystyle\int\limits_{\pi/4}^{\pi/2}}
\left[  \frac{\pi}{4}-2\arctan\left(  \frac{-\sin(\varphi)+\sqrt{-\cos
(\varphi)^{2}+2\sin(\varphi)^{2}}}{\cos(\varphi)+\sin(\varphi)}\right)
\right]  d\varphi\\
&  =0.2815898507...
\end{align*}
where $\lambda=\arccos(\sqrt{2/3})$ and, conditional on this, the obtuseness
probability is%
\[
1-\frac{8}{\Delta\pi^{2}}%
{\displaystyle\int\limits_{\pi/4}^{\pi/2}}
\left[  \frac{\pi}{4}-2\arccos\left(  \frac{1}{\sqrt{2}\sin(\varphi)}\right)
\right]  d\varphi=0.6597451305....
\]

\section{Integration Details}

Starting with the bivariate density for sides $a$, $c$ in Section 5, let
$x=2a^{2}-1$, then
\[%
\begin{array}
[c]{lll}%
dx=4a\,da, &  & 1-a^{2}=1-\left(  \dfrac{1+x}{2}\right)  =\dfrac{1-x}{2}%
\end{array}
\]
hence%
\[%
\begin{array}
[c]{lll}%
\dfrac{\sqrt{2}}{4\sqrt{1+x}}dx=da, &  & 4a^{2}\left(  1-a^{2}\right)
=2(1+x)\left(  \dfrac{1-x}{2}\right)  =1-x^{2}%
\end{array}
\]
and thus the density becomes%
\[
\dfrac{4}{\pi^{2}}\frac{\sqrt{2}}{\sqrt{1-x}}\dfrac{c}{\sqrt{\left(
1-x^{2}\right)  -\left(  1-c^{2}\right)  ^{2}}}\dfrac{\sqrt{2}}{4\sqrt{1+x}%
}=\dfrac{2}{\pi^{2}}\dfrac{c}{\sqrt{1-x^{2}}\sqrt{\left[  1-\left(
1-c^{2}\right)  ^{2}\right]  -x^{2}}}.
\]
Integrating with respect to $x$, formula 3.152.7 in \cite{GR-UnifEqua} can be
applied to obtain the univariate density for $c$. \ To compute
$\operatorname*{E}(a\,c)$, examine instead \
\[
\dfrac{2}{\pi^{2}}\dfrac{a\,c^{2}}{\sqrt{1-x^{2}}\sqrt{\left[  1-\left(
1-c^{2}\right)  ^{2}\right]  -x^{2}}}=\dfrac{\sqrt{2}}{\pi^{2}}\dfrac{c^{2}%
}{\sqrt{1-x}\sqrt{\left[  1-\left(  1-c^{2}\right)  ^{2}\right]  -x^{2}}}
\]
and apply formula 3.131.3 in \cite{GR-UnifEqua}. The analogous expression for
$\operatorname*{E}(a\,c)$ in Section 4:%
\[
\dfrac{4}{\pi}\dfrac{a^{2}c^{2}}{\sqrt{\left(  1-x^{2}\right)  -\left(
1-c^{2}\right)  ^{2}}}\dfrac{\sqrt{2}}{4\sqrt{1+x}}=\dfrac{1}{\sqrt{2}\pi
}\dfrac{c^{2}\sqrt{1+x}}{\sqrt{\left[  1-\left(  1-c^{2}\right)  ^{2}\right]
-x^{2}}}
\]
can be integrated via formula 3.141.17 in \cite{GR-UnifEqua}.

Starting with the bivariate density for sides $a$, $b$ in Section 6, let%
\[%
\begin{array}
[c]{ccc}%
x=\sqrt{a^{2}+b^{2}}, &  & y=a/b
\end{array}
\]
then%
\[%
\begin{array}
[c]{ccccc}%
\left\vert J\right\vert =\dfrac{1+y^{2}}{x}, &  & a=\dfrac{x\,y}{\sqrt
{1+y^{2}}}, &  & b=\dfrac{x}{\sqrt{1+y^{2}}}%
\end{array}
\]
and the density becomes separable:%
\[
\frac{1}{\sqrt{1-x^{2}}\left(  1+y^{2}\right)  }
\]
although the region of integration initially appears unmanageable. \ The
inequalities%
\[
\dfrac{x\left\vert y-1\right\vert }{\sqrt{1+y^{2}}}<\sqrt{1-x^{2}}%
<\dfrac{x\left(  y+1\right)  }{\sqrt{1+y^{2}}}
\]
become%
\[
\frac{\sqrt{1+y^{2}}}{\sqrt{2}\sqrt{1+y+y^{2}}}<x<\frac{\sqrt{1+y^{2}}}%
{\sqrt{2}\sqrt{1-y+y^{2}}}
\]
and, integrating with respect to $x$, we obtain%
\begin{align*}
&  \frac{1}{1+y^{2}}\left[  \arcsin\left(  \frac{\sqrt{1+y^{2}}}{\sqrt{2}%
\sqrt{1-y+y^{2}}}\right)  -\arcsin\left(  \frac{\sqrt{1+y^{2}}}{\sqrt{2}%
\sqrt{1+y+y^{2}}}\right)  \right] \\
&  =\frac{1}{1+y^{2}}\left[  \arctan\left(  \frac{\sqrt{1+y^{2}}}{\left\vert
1-y\right\vert }\right)  -\arctan\left(  \frac{\sqrt{1+y^{2}}}{1+y}\right)
\right] \\
&  =\left\{
\begin{array}
[c]{lll}%
\dfrac{1}{1+y^{2}}\arctan\left(  y\sqrt{1+y^{2}}\right)  &  & \text{if
}0<y<1,\\
\dfrac{1}{1+y^{2}}\arctan\left(  \dfrac{\sqrt{1+y^{2}}}{y^{2}}\right)  &  &
\text{if }1<y<\infty
\end{array}
\right.
\end{align*}
since%
\[
2\left(  1\pm y+y^{2}\right)  -\left(  1+y^{2}\right)  =1\pm2y+y^{2}%
=\left\vert 1\pm y\right\vert ^{2}
\]
and by the addition formula for the arctangent. \ Symmetry between $0<y<1$ and
$1<y<\infty$ allows us to focus on the former when integrating over $y$ and to
multiply the final result by two.

\section{Addendum}

I\ am grateful to M. Larry Glasser for reducing $\operatorname*{E}(c)$ in
Section 5 to%
\[
\frac{4\sqrt{2}}{3\pi}\,_{3}F_{2}\left(  \dfrac{1}{2},\dfrac{1}{2},\dfrac
{3}{2};\dfrac{5}{4},\dfrac{7}{4};1\right)
\]
where $_{3}F_{2}$ is a hypergeometric function
\[
_{3}F_{2}\left(  \dfrac{1}{2},\dfrac{1}{2},\dfrac{3}{2};\dfrac{5}{4},\dfrac
{7}{4};x\right)  =\frac{3}{4\sqrt{2\pi}}%
{\displaystyle\sum\limits_{n=0}^{\infty}}
\frac{\Gamma(n+1/2)^{2}\Gamma(n+3/2)}{\Gamma(n+5/4)\Gamma(n+7/4)}\frac{x^{n}%
}{n!}%
\]
and Michael S. Milgram \cite{Mg-UnifEqua} for pointing out the further
simplification%
\[
\frac{\pi^{4}+8\Gamma\left(  \frac{5}{8}\right)  ^{4}\Gamma\left(  \frac{7}%
{8}\right)  ^{4}}{2\pi^{3}\Gamma\left(  \frac{5}{8}\right)  ^{2}\Gamma\left(
\frac{7}{8}\right)  ^{2}}%
\]
which follows from equating $_{3}F_{2}(1/2,1/2,3/2;5/4,7/4;1)$ with
$\Omega_{1,1}(1/2,1/2,5/4)$ in Example 13 of \cite{Ch-UnifEqua}.  Much more
relevant material can be found at \cite{Fi3-UnifEqua}, including experimental
computer runs that aided theoretical discussion here.

\end{document}